\documentclass{article}
\usepackage{titling}
\usepackage{datetime}
\usepackage{float}
\usepackage{multirow}
\usepackage{makecell}
\usepackage{array}

\usepackage{xcolor}
\usepackage{hyperref}
\hypersetup{
    colorlinks,
    linkcolor={red!50!black},
    citecolor={blue!50!black},
    urlcolor={blue!80!black}
}

\usepackage{enumitem}
\makeatletter
\def\namedlabel#1#2{\begingroup
    #2%
    \def\@currentlabel{#2}%
    \phantomsection\label{#1}\endgroup
}
\makeatother

\usepackage[top=2.5cm,bottom=2.5cm,right=2cm,left=2cm]{geometry}

\usepackage{amsmath, amssymb, mathtools, amsthm}
\usepackage{empheq}
\usepackage[thinc]{esdiff}
\usepackage{accents}

\DeclareMathOperator*{\R}{\mathbb{R}}

\DeclareMathOperator*{\Natural}{\mathbb{N}}

\usepackage{tikz}
\usepackage{bbold}
\usepackage[font=small]{caption}
\usepackage[font=small]{subcaption}
\usepackage[section]{placeins}

\usepackage{algpseudocode}
\usepackage{algorithm}

\allowdisplaybreaks

\usepackage{mdframed}
\newmdenv[
  topline=false,
  bottomline=false,
  rightline=false,
  skipabove=\topsep,
  skipbelow=\topsep,
  leftmargin=2pt,
  rightmargin=0pt,
  innertopmargin=0pt,
  innerbottommargin=0pt
]{sideline}

\usepackage[version=4]{mhchem}
\def\arrvline{\hfil\kern\arraycolsep\vline\kern-\arraycolsep\hfilneg}

\usepackage[sorting=none, style=vancouver]{biblatex}
\bibliography{main}

\AtBeginEnvironment{proof}{\footnotesize \linespread{0.7}}

\title{The generator gradient estimator is an adjoint state method for stochastic differential equations}
\author{Quentin Badolle$^1$, Ankit Gupta$^1$, Mustafa Khammash$^{1*}$}
\date{\monthname{} \the\year}
\date{
$^1$Department of Biosystems Science and Engineering, ETH Zurich, Basel, Switzerland \\
$^*$ mustafa.khammash@bsse.ethz.ch
}

\begin{document}

\maketitle

\begin{abstract}
Motivated by the increasing popularity of overparameterized Stochastic Differential Equations (SDEs) like Neural SDEs, Wang, Blanchet and Glynn recently introduced the \emph{generator gradient estimator}, a novel unbiased stochastic gradient estimator for SDEs whose computation time remains stable in the number of parameters~\cite{wang2024efficient}. In this note, we demonstrate that this estimator is in fact an adjoint state method, an approach which is known to scale with the number of states and \emph{not} the number of parameters in the case of Ordinary Differential Equations~(ODEs). In addition, we show that the generator gradient estimator is a close analogue to the exact Integral Path Algorithm~(eIPA) estimator which was introduced by Gupta, Rathinam and Khammash for a class of Continuous-Time Markov Chains (CTMCs) known as stochastic chemical reactions networks (CRNs)~\cite{gupta2018estimation}.
\end{abstract}

Given a terminal time $T$, we consider a family of diffusions $\{X_{\theta}^{x}(t, s) \in \mathbb{R}^{d} : s \in [t, T]\}$ indexed by an initial condition $x \in \mathbb{R}^{d}$ at time $t$ and a parameter $\theta \in \Theta \subset \mathbb{R}^n$. The dynamics of the process $(X_{\theta}^{x}(t, s))_{s\in [t, T]}$ are generated by an Itô SDE given by:

\begin{equation}
\label{eq:dynamics_sde}
X_{\theta}^{x}(t,s) = x + \int_{t}^{s} \mu_{\theta}(r, X_{\theta}^{x}(t,r))dr + \int_{t}^{s}\sigma_{\theta}(r, X_{\theta}^{x}(t,r))dB(r),
\end{equation}

where $\mu_{\theta}$ is the drift term and $\sigma_{\theta}$ the volatility term. We are interested in estimating the gradient of:

\begin{equation}
\label{eq:objective_sde}
v_{\theta}(t,x) \coloneqq E\Bigg[\int_{t}^{T}\rho_{\theta}(s, X_{\theta}^{x}(t,s))ds + g_{\theta}(X_{\theta}^{x}(t,T))\Bigg],
\end{equation}

where $\rho_{\theta}$ represents the reward rate and and $g_{\theta}$ the terminal reward.

\section{The generator gradient estimator as an adjoint state method}

As in~\cite{wang2024efficient}, without loss of generality, we focus on the gradient of $v_{\theta}$ at time $t=0$. For ease of exposition, we introduce $X_{\theta}^{x}(t) \coloneq X_{\theta}^{x}(0, t)$. We denote $a_{\theta}$ the diffusion matrix given by:

\begin{equation}
\label{eq:diffusion_matrix}
a_{\theta}(t, x) \coloneqq \frac{1}{2}\sigma_{\theta}(t,x)\sigma_{\theta}(t,x)^{\intercal},
\end{equation}

which we use to define the generator $\mathcal{L}_{\theta}$ of $(X_{\theta}^{x}(t))$ as:

\begin{equation}
\label{eq:generator_sde}
\mathcal{L}_{\theta}f(t, x) \coloneqq (\nabla f(t,x))^{\intercal}\mu_{\theta}(t,x) + \text{Tr}(\nabla^{2}f(t,x)a_{\theta}(t,x)),
\end{equation}

where $f$ is a twice differentiable function. $\nabla$ corresponds to a space gradient, $\nabla^2$ to a space hessian and $\text{Tr}$ to the trace operator.  It was rigorously shown in~\cite{wang2024efficient} that the derivative of $v_{\theta}$ can be expressed as:

\begin{equation}
\label{eq:derivative_sde}
\partial_{\theta_{i}}v_{\theta}(0,x) = E\Bigg[\int_{0}^{T} \Big(\partial_{\theta_{i}}\mathcal{L}_{\theta}v_{\theta}(t, X_{\theta}^{x}(t))+\partial_{\theta_{i}}\rho_{\theta}(t, X_{\theta}^{x}(t))\Big)dt+\partial_{\theta_{i}} g_{\theta}(X_{\theta}^{x}(T))\Bigg],
\end{equation}

where $\partial_{\theta_{i}} \mathcal{L}_{\theta}$ has been defined as:

\begin{equation}
\label{eq:generator_derivative_sde}
\partial_{\theta_{i}}\mathcal{L}_{\theta}f(t, x) \coloneqq (\nabla f(t,x))^{\intercal}\partial_{\theta_{i}}\mu_{\theta}(t,x) + \text{Tr}(\nabla^{2}f(t,x)\partial_{\theta_{i}}a_{\theta}(t,x)),
\end{equation}

and $\partial_{\theta_{i}}$ denotes the element-wise partial derivative with respect to the $i$-th coordinate of $\theta$.\newline

We first provide an informal proof of eq.~\eqref{eq:derivative_sde} to establish the connection with the adjoint state method (see~\cite{cea1986conception, plessix2006review} for a presentation of this approach). Similarly to~\cite{wang2024efficient}, we start from the Feynman-Kac formula which states that, under sufficient regularity conditions, $v_{\theta}$ is the solution of the Partial Differential Equation (PDE) given by (see section 7 in chapter 5 of~\cite{karatzas2014brownian}):

\begin{equation}
\label{eq:fc_pde_objective}
\partial_{t} v_{\theta} + \mathcal{L}_{\theta} v_{\theta} + \rho_{\theta} = 0, \quad v_{\theta}(T, \cdot) = g_{\theta},
\end{equation}

where $\partial_{t}$ denotes a time derivative. Assuming enough smoothness, we formally differentiate the PDE in eq.~\eqref{eq:fc_pde_objective} with respect to $\theta_{i}$ to obtain:

\begin{equation}
\label{eq:fc_pde_objective_derivative}
\partial_{t}\partial_{\theta_{i}} v_{\theta} + \mathcal{L}_{\theta} \partial_{\theta_{i}}v_{\theta} + \partial_{\theta_{i}}\mathcal{L}_{\theta}  v_{\theta}+ \partial_{\theta_{i}} \rho_{\theta} = 0, \quad \partial_{\theta_{i}} v_{\theta}(T, \cdot) = \partial_{\theta_{i}} g_{\theta}.
\end{equation}

Given two functions $f$ and $g$ defined on $\mathbb{R}^d$ and taking values in $\mathbb{R}$, we define the inner product $\langle \cdot, \cdot \rangle$ as:

$$
\langle f, g \rangle \coloneqq \int_{\mathbb{R}^d} f(x) g(x)dx.
$$

Introduce a family of probability mass functions $\{p(t, \cdot) : t \in [0, T]\}$ indexed by time, which we leave unspecified for the moment. Taking the inner product of eq.~\eqref{eq:fc_pde_objective_derivative} with $p(t, \cdot)$ and integrating over $[0, T]$ leads to:

\begin{align}
\label{eq:int_fc_pde_objective_derivative}
\begin{split}
\int_{0}^{T}\Big\langle\partial_{t}\partial_{\theta_{i}} v_{\theta}(t, \cdot),p(t, \cdot)\Big\rangle dt &+ \int_{0}^{T}\Big\langle\mathcal{L}_{\theta} \partial_{\theta_{i}}v_{\theta}(t, \cdot),p(t, \cdot)\Big\rangle dt\\
&+ \int_{0}^{T}\Big\langle\partial_{\theta_{i}}\mathcal{L}_{\theta} v_{\theta}(t, \cdot),p(t, \cdot)\Big\rangle dt+ \int_{0}^{T}\Big\langle\partial_{\theta_{i}} \rho_{\theta}(t, \cdot),p(t, \cdot)\Big\rangle dt = 0.
\end{split}
\end{align}

By integration by parts, the first term on the left-hand side of eq.~\eqref{eq:int_fc_pde_objective_derivative} can be rewritten as:

\begin{align}
\int_{0}^{T}\Big\langle\partial_{t}\partial_{\theta_{i}} v_{\theta}(t, \cdot),p(t, \cdot)\Big\rangle dt
&= \int_{\mathbb{R}^d}\int_{0}^{T}\partial_{t}\partial_{\theta_{i}} v_{\theta}(t, x')p(t, x')dtdx' \nonumber\\
&= \int_{\mathbb{R}^d}\bigg[\partial_{\theta_{i}} v_{\theta}(t, x')p(t, x')\bigg]_{0}^{T}dx' - \int_{\mathbb{R}^d}\int_{0}^{T}\partial_{\theta_{i}} v_{\theta}(t, x')\partial_{t}p(t, x')dtdx' \nonumber\\
&= \int_{\mathbb{R}^d}\Big(\partial_{\theta_{i}} v_{\theta}(T, x')p(T, x') - \partial_{\theta_{i}} v_{\theta}(0, x')p(0, x')\Big)dx' - \int_{0}^{T}\Big\langle\partial_{\theta_{i}} v_{\theta}(t, \cdot),\partial_{t} p(t, \cdot)\Big\rangle dt. \label{eq:intermediate_step}\\
\intertext{Now choose $p(t,\cdot)$ to be $p_{\theta}(t,\cdot)$, the probability mass function of $(X_{\theta}^{x}(t))$ at time $t$ and write $\delta_x$ the Kronecker delta. Observe in eq.~\eqref{eq:intermediate_step} that:}
\noalign{\centering $\partial_{\theta_{i}} v_{\theta}(T, \cdot) = \partial_{\theta_{i}} g_{\theta}$, \quad $\partial_{\theta_{i}} v_{\theta}(0, \cdot)p_{\theta}(0, \cdot) = \partial_{\theta_{i}}v_{\theta}(0, \cdot) \delta_x$.}
\intertext{Using the Fokker-Planck equation, we then have:}
\int_{0}^{T}\Big\langle \partial_{t}\partial_{\theta_{i}}v_{\theta}(t, \cdot),p_{\theta}(t, \cdot)\Big\rangle dt 
&= \Big\langle \partial_{\theta_{i}}g_{\theta},p_{\theta}(T, \cdot)\Big\rangle -\partial_{\theta_{i}}v_{\theta}(0, x) -
\int_{0}^{T}\Big\langle\partial_{\theta_{i}} v_{\theta}(t, \cdot),\mathcal{L}_{\theta}^{*}p_{\theta}(t, \cdot)\Big\rangle dt\nonumber\\
&= \Big\langle \partial_{\theta_{i}} g_{\theta},p_{\theta}(T, \cdot)\Big\rangle -\partial_{\theta_{i}}v_{\theta}(0, x) -
\int_{0}^{T}\Big\langle\mathcal{L}_{\theta}\partial_{\theta_{i}} v_{\theta}(t, \cdot),p_{\theta}(t, \cdot)\Big\rangle dt\label{eq:first_term_int_fc},
\end{align}

where $\mathcal{L}_{\theta}^{*}$ is the adjoint operator of $\mathcal{L}_{\theta}$. Notice that the third term in eq.~\eqref{eq:first_term_int_fc} is the opposite of the second term in eq.~\eqref{eq:int_fc_pde_objective_derivative}. Therefore, plugging eq.~\eqref{eq:first_term_int_fc} in eq.~\eqref{eq:int_fc_pde_objective_derivative}, we get:

\begin{align}
\label{eq:derivative_sde_scalar_product}
\partial_{\theta_{i}}v_{\theta}(0, x) 
&= \int_{0}^{T}\Big\langle\partial_{\theta_{i}}\mathcal{L}_{\theta} v_{\theta}(t, \cdot),p_{\theta}(t, \cdot)\Big\rangle dt+ \int_{0}^{T}\Big\langle\partial_{\theta_{i}} \rho_{\theta}(t, \cdot),p_{\theta}(t, \cdot)\Big\rangle dt + \Big\langle \partial_{\theta_{i}} g_{\theta},p_{\theta}(T, \cdot)\Big\rangle,
\end{align}

where $p_{\theta}(t,\cdot)$ is the solution of an adjoint state equation which is here the Fokker-Planck equation with initial condition $\delta_{x}$. Eq.~\eqref{eq:derivative_sde_scalar_product} is just another way to write eq.~\eqref{eq:derivative_sde}. This informal derivation draws an explicit connection between eq.~\eqref{eq:derivative_sde} used for the  generator gradient estimator and the well-known adjoint state method.\newline

The generator gradient estimator based on eq.~\eqref{eq:derivative_sde} relies on estimates of $\nabla v_{\theta}$ and $\nabla^2   v_{\theta}$ obtained from auxiliary pathwise differentiation estimators. Remarkably, the number of such auxiliary estimators scales with the number of states $d$ and \emph{not} the number of parameters $n$, exactly as with the adjoint state method in the case of ODEs.

\section{The generator gradient estimator as an analogue to the exact Integral Path Algorithm~(eIPA) estimator}

Observe that the derivation of eq.~\eqref{eq:derivative_sde} does not rely on the explicit expression for the generator $\mathcal{L}_{\theta}$. In fact, a similar expression was obtained in the context of a class of CTMCs known as CRNs~\cite{gupta2018estimation}~(see~\cite{erdi1989mathematical, anderson2015stochastic, wilkinson2018stochastic} for an introduction to these models). To make the parallel explicit, we need to introduce some notations. Let us consider a network with $d$ molecular species. The state of the system at any time can be described by a vector in $\Natural^d$ whose $i$-th component corresponds  to the number of molecules of the $i$-th species. The chemical species interact through $m$ chemical reactions and every time the $k$-th reaction fires, the state of the system is displaced by the $d$-dimensional stoichiometric vector $\zeta_{k}\in\mathbb{Z}^d$. The time-homogeneous propensity function $\lambda_{\theta} = (\lambda_{\theta, k})_{k\in [\![1, m]\!]}$ parameterised by a parameter $\theta \in \Theta \subset \R^n$ depends on the state of the system $x \in \Natural^d$.  Given a terminal time $T$, we consider a family of Markov jump processes  $\{X_{\theta}^{x}(t) : t \in [0, T]\}$ indexed by an initial condition $x$ at time $t=0$ and the parameter $\theta$. The generator of $(X_{\theta}^{x}(t))$ is given by:

\begin{equation}
\label{eq:generator_ctmc}
\mathcal{L}_{\theta}f(x) \coloneqq \sum_{k=1}^{m} \lambda_{\theta,k}(x)\Delta_{\zeta_k}f(x),
\end{equation}

where, given a vector $z \in \R^{d}$, $\Delta_{z}$ is the spatial finite difference operator: $\Delta_{z} f(x) \coloneqq f(x + z) - f(x)$. Given a collection of independent, unit-rate Poisson processes $\{(Y_k(t))\}_{k \in [\![1,m]\!]}$, we associate to each reaction $k$ a counting process $(R_{\theta, k}(t))$ defined as:

\begin{equation}
\label{eq:reaction_count}
R_{\theta, k}(t) \coloneqq Y_{k}\bigg(\int_0^t \lambda_{\theta, k}(X_{\theta}^{x}(s)) ds\bigg).  
\end{equation}

The random time change representation of $(X_{\theta}^{x}(t))$ is given by~\cite{anderson2015stochastic, kurtz1980representations}:

\begin{equation}
\label{eq:dynamics_ctmc}    
X_{\theta}^{x}(t) = x + \sum_{k=1}^{m} \zeta_{k} R_{\theta, k}(t).
\end{equation}

We redefine $v_{\theta}$ as:

\begin{equation}
\label{eq:objective_ctmc}
v_{\theta}(t,x) \coloneqq E\Big[g(X_{\theta}^{x}(T-t))\Big],
\end{equation}



where $g$ corresponds to a terminal reward. Theorem 3.3 in~\cite{gupta2018estimation} gives an expression for $\partial_{\theta_{i}}v_{\theta}$ which is rigorously derived. Using the time-homogeneity of $(X_{\theta}^{x}(t))$ to rewrite the formula given there, we get:

\begin{equation}
\label{eq:derivative_ctmc}
\partial_{\theta_{i}}v_{\theta}(0,x) = \sum_{k=1}^{m} E\Bigg[\int_{0}^{T} \Delta_{\zeta_k}v_{\theta}(t, X_{\theta}^{x}(t)) \partial_{\theta_{i}}\lambda_{\theta, k}(X_{\theta}^{x}(t))dt\Bigg].
\end{equation}

This expression can straightforwardly be extended to include some dependence of the terminal reward $g$ on $\theta$, a reward rate $\rho_{\theta}$ and time-dependent propensities. To make the comparison between eq.~\eqref{eq:derivative_ctmc} and eq.~\eqref{eq:derivative_sde}-\eqref{eq:generator_derivative_sde} transparent, let us rewrite eq.~\eqref{eq:derivative_sde} for a time-homogeneous SDE, in the absence of reward rate and in the case when the terminal reward is independent of $\theta$:

\begin{equation}
\label{eq:derivative_sde_simplified}
\partial_{\theta_{i}}v_{\theta}(0,x) = E\Bigg[\int_{0}^{T} \Big((\nabla v_{\theta}(t, X_{\theta}^{x}(t)))^{\intercal}\partial_{\theta_{i}}\mu_{\theta}(X_{\theta}^{x}(t)) + \text{Tr}(\nabla^{2}v_{\theta}(t,X_{\theta}^{x}(t))\partial_{\theta_{i}}a_{\theta}(X_{\theta}^{x}(t))\Big)dt\Bigg].
\end{equation}

In particular, the spatial finite difference terms $\Delta_{\zeta_k}v_{\theta}$ in eq.~\eqref{eq:derivative_ctmc} are direct analogues to $\nabla v_{\theta}$ and $\nabla^{2} v_{\theta}$ in eq.~\eqref{eq:derivative_sde_simplified}, and the parametric partial derivatives $\partial_{\theta_{i}} \lambda_{\theta, k}$ are the equivalent of $\partial_{\theta_{i}}\mu_{\theta}$ and $\partial_{\theta_{i}}a_{\theta}$.\newline

Estimation of $\partial_{\theta_{i}}v_{\theta}$ based on eq.~\eqref{eq:derivative_ctmc} as part of the exact Integral Path Algorithm~(eIPA) introduced in~\cite{gupta2018estimation} shares key ideas with the generator gradient estimator. First, eIPA replaces the unknown quantities $v_{\theta}$ by so called auxiliary processes which should be compared to the estimates of $\partial_{\theta_{i}}v_{\theta}$ obtained from pathwise differentiation estimators~(see the discussion around eq.~(2.6) in~\cite{wang2024efficient} and subsection 3.4 in~\cite{gupta2018estimation}). Secondly, to control the computation time per sample, the auxiliary processes of eIPA are generated at a given jump time for reaction $k$ only if a certain Bernoulli random variable equals $1$~(again, see subsection 3.4 in~\cite{gupta2018estimation}). This is similar in spirit to the integral randomisation strategy used by the generator gradient estimator~(see the paragraph related to eq.~(2.9) in~\cite{wang2024efficient}). Strikingly, it was illustrated in~\cite{gupta2018estimation} that the eIPA estimator exhibits low variance when compared to other unbiased estimators, as is the case for the generator gradient estimator~(see the numerical examples in section 4 of~\cite{wang2024efficient} and~\cite{gupta2018estimation}).

\printbibliography

\end{document}